\documentclass[11pt,twoside,draft]{article}
\usepackage{mathbbold}

\usepackage{amsmath, amssymb, mathrsfs, graphicx, mathabx}

\def\az{\alpha}  \def\bz{\beta}
    \def\dz{\delta}
    \def\fz{\varphi}
\def\gz{\gamma}  
\def\lz{\lambda}

        \def\uz{\theta}

 \def\ddz{\Delta}
    
  \def\ooz{\Omega}
\def\ppz{\Pi}

\def\qqd{\qquad}

\newcommand{\mathsym}[1]{{}}

\def\le{\leqslant}
\def\ge{\geqslant}

\font\cms=cmss9 scaled \magstep1

\def\nnd{\noindent}

\def\thm{\nnd\bg{thm1}}
\def\crl{\nnd\bg{crl1}}
\def\lmm{\nnd\bg{lmm1}}

\def\xmp{\nnd\bg{xmp1}}

\def\dethm{\end{thm1}}
\def\decrl{\end{crl1}}
\def\delmm{\end{lmm1}}

\def\dexmp{\end{xmp1}}

\def\prf{\medskip \noindent {\bf Proof}. }
\def\qed{\text{\quad $\Box$}}
\def\deprf{\qed\medskip}
\def\bg{\begin}
\def\be{\bg{equation}}
\def\de{\end{equation}}

\def\dear{\end{eqnarray}}
\def\lb{\label}
\def\ct{\cite}

\newcommand{\rf}[2]{[\ref{#1}; #2]}

\def\den{\end{enumerate}}

\def\d{\text{\rm d}}

\allowdisplaybreaks[4]

\begin{document}


\thispagestyle{empty}
\renewcommand{\thefootnote}{\fnsymbol{footnote}}

\noindent {Chinese Journal of Applied Probability and Statistics 2015}\newline

\begin{center}
{\bf\Large Practical Criterion for Uniqueness of $Q$-processes}
\vskip.15in {Mu-Fa Chen}
\end{center}
\begin{center} (Beijing Normal University)\\
\vskip.1in December 31, 2014
\end{center}
\vskip.1in

\markboth{\sc Mu-Fa Chen}{\sc Practical Criterion for Uniqueness}



\date{}



\begin{abstract}
The note begins with a short story on seeking for a practical sufficiency
theorem for the uniqueness of time-continuous Markov jump processes, starting around 1977.
The general result was obtained in 1985 for the processes with general
state spaces. To see the sufficient conditions are sharp, a dual criterion
for non-uniqueness was obtained in 1991. This note is restricted however
to the discrete state space (then the processes are called $Q$-processes or Markov chains), for which the sufficient conditions just mentioned are showing at the end of the note to be necessary. Some examples are included to illustrate that
the sufficient conditions either for uniqueness or for non-uniqueness are not only powerful but also sharp.
\end{abstract}

\nnd {\small 2000 {\it Mathematics Subject Classification}: 60J27}

\nnd {\small {\it Key words and phrases}. Criterion; uniqueness; Markov chain;
Markov jump process.}

\bigskip

Let $E$ be a countable set with elements $i, j, k, \cdots$. A matrix $Q=(q_{ij}: i, j\in E)$ is called a
$Q$-matrix if its non-diagonals are nonnegative and $\sum_{j\in E}q_{ij}\le 0$ for every $i\in E$.
Throughout this note, we restrict ourselves to the special case that the $Q$-matrix is totally stable
$q_i:=-q_{ii}<\infty$ and conservative $q_i=\sum_{j\ne i}q_{ij}$ for every $i\in E$. It is called bounded if $\sup_{i\in E}q_i<\infty$. For a given $Q$-matrix $Q=(q_{ij})$ on $E$,
a sub-Markovian semigroup $\{P(t)=(p_{ij}(t): i, j\in E)\}_{t\ge 0}$ is
called a $Q$-process if
$$\frac{\d}{\d t} P(t)\bigg|_{t=0}=Q\qquad\text{\rm (pointwise)}.$$
The $Q$-processes may not be unique
in general, but there always exists the minimal one, due to Feller (1940)\,\rf{wf40}{Theorem 1},
denoted by $P^{\min}(t)=\big(p_{ij}^{\min}(t): i, j\in E\big)$.
For more than half-century ago, some criteria for the uniqueness were
known.

\thm{\cms
The $Q$-process is unique (equivalently, the minimal process $P^{\min}(t)$ is not explosive) iff one of the following equivalent conditions holds:
\begin{itemize} \setlength{\itemsep}{0cm}%
\item[(C1)] $\sum_{j\in E} p_{ij}^{\min}(t)=1$ for every $i\in E$ and $t\ge 0$.
\item[(C2)] $\sum_{n=1}^\infty q_{X^{\min}(\tau_n)}^{-1}=\infty$, ${\mathbb P}_i$-a.s., where $\tau_n$ is the
  $n$th jump time of the minimal process $\{X^{\min}(t): t\ge 0\}$ corresponding to   $P^{\min}(t)$.
\item[(C3)] The equation
\be (\lambda I-Q) u=0, \qquad 0\le u\le 1,\lb{01}\de
has only zero solution for some (equivalently, for all) $\lambda>0$.
\end{itemize}
}\dethm
Criterion (C1) goes back to Feller (1940)\,\cite{wf40}. Criterion (C2)
is due to Dobrushin (1952)\,\cite{drl52}. Criterion (C3) is due to
Feller (1957)\,\cite{wf57} and Reuter (1957) \cite{gehr57}.
Refer also to Chung (1960)\,\rf{ckl67}{Part II, \S 19, Theorem 1}, or \rf{gs75}{Chap. 3, \S 2,
Theorems 3 and 4}.

The earlier Criterion (C1) often requires a further effort in practice, rather than
a direct application. In particular, the proof of the powerful sufficiency theorem (Theorem \ref{t1} below) is based on it.

 Criterion (C2) is effective in some cases. For instance in the
simplest case that $M:=\sup_{i\in E}q_i<\infty$, since
$$\sum_{n=1}^\infty q_{X^{\min}(\tau_n)}^{-1}\ge \sum_{n=1}^\infty M^{-1}=\infty,$$
we obtain the uniqueness of the processes. For pure birth process (i.e.,
$q_{i, i+1}>0$ and $q_{ij}=0$ for all $j\ne i$, $i, j\ge 0$), Criterion (C2)
says that the process is unique iff
\be \sum_{n=1}^\infty \frac{1}{q_{n, n+1}}=\infty. \lb{02}\de
Besides, if the minimal process is recurrent, then the term $q_k^{-1}$ will appears
infinitely often in the summation, hence the process should be unique according to
the criterion.

Criterion (C3) is more effective once equation (\ref{01}) is solvable.
More precisely, it is the case if the exit boundary consists at most a single point,
for instance the pure birth processes,
the birth--death processes or more general the single birth processes
(i.e., for $j>i\ge 0$, $q_{ij}>0$ iff $j=i+1$; for $0\le j<i$, $q_{ij}$
is nonnegative but free). We will come back this story soon.

However, the next model stopped our study for several years at the beginning
of the study (1977-1978) on non-equilibrium particle systems. To state our model,
we use operator $\ooz$ instead of the matrix $Q$:
$$\ooz f(i)=\sum_{j\in E}q_{ij}(f_j-f_i),\qqd i\in E.$$
Of course, in this case, $\ooz f=Qf$. For a Markov chain on a countable
set $E$, by a transform, one often assumes that $E$ is simply the set
${\mathbb Z}_+=\{0, 1, \cdots\}$. However, such a transform ignores the
original geometry of $E$ and may not be convenient in multidimensional case.
To state our model, we need some notation. Let $i=(i_u: u\in S)$ and define
its updates $i^{u\pm}$ and $i^{u, v}$ as follows:
$${i_w^{u\pm}=\bg{cases} i_u\pm 1\quad &w=u\\ i_w &w\ne u,\end{cases}}\qquad
{i_w^{u, v}=\bg{cases} i_u-1 \quad &w=u\\ i_v+1 & w=v\\ i_w & w\ne u, v,\end{cases}}
\qquad w\in S.$$

\xmp\lb{t0}{\cms[Schl\"ogl's second model]}\;
{\cms Let $S$ be a finite set and $E={\mathbb Z}_+^S$. Define a Markov chain on
$E$ with operator
$$\aligned
\ooz f(i)&= \sum_{u\in S}\Big\{b(i_u)\big[
f\big(i^{u+}\big)-f(i)\big] +a(i_u)\big[
f\big(i^{u-}\big)-f(i)\big]\Big\}\\
&\quad+\!\sum_{u, v} i_u p(u, v)
\big[f\big(i^{u, v}\big)-f(i)\big],\qquad  i=(i_u: u\in S)\in E,
\endaligned$$
where $(p(u, v): u, v\in S)$ is a ``simple'' random walk on $S$, and
$$\aligned
b(k)&=\bz_0+\bz_2k (k-1),\quad \bz_0, \bz_2>0,\\
a(k)&=\dz_1 k + \dz_3 k (k-1)(k-2), \quad \dz_1, \dz_3 >0.
\endaligned$$
}\dexmp

Here in the first sum of $\ooz$, in each vessel $u$, there is a birth--death
process with birth rate $b(k)$ and death rate $a(k)$, respectively. This is called the reaction part of the model. The reactions in different vessels are independent. In the second sum of $\ooz$, a particle from vessel $u$ moves to vessel $v$. This is called the diffusion part of the model. Thus, it is actually a finite-dimensional reaction--diffusion processes. Replacing the finite $S$ with $S={\mathbb Z}^d$, we obtain formally an operator of infinite-dimensional reaction--diffusion process
which is a typical model from the non-equilibrium statistical physics.
Even though the large systems are quite popular today, in that period, it was
rather unusual to study such a non-equilibrium system. Our original program is
to rebuild the mathematical ground of non-equilibrium statistical physics (cf. \rf{cmf04}{Part IV}. An earlier paper on this topic appeared in 1985 \cite{cmf85}). For this,
the model is meaningful only if it is ergodic in every finite dimension.
Thus, the finite dimensional model consists the first doorsill of our program.

In 1983, the author and Yan \cite{ysjcmf86}, using a comparison technique, overcame this doorsill, based on a systemic study on the single birth processes.
To which, we obtained explicit criteria not only for uniqueness but
also for ergodicity and so on. This goes back to \cite{ysjcmf86, cmf86a}.
Refer to \ct{cmf04} for updates and to \cite{chzh} for a unified treatment.
After two more years, using an approximating approach, we obtained a powerful sufficiency theorem as stated below.

\thm\lb{t1}{\cms[Uniqueness criterion]}\;
{\cms Let $Q=(q_{ij})$ be a $Q$-matrix on a coun\-table set $E$. Then the corresponding $Q$-process is unique iff the following two conditions
hold simultaneously.
\begin{itemize}\setlength{\itemsep}{0cm}
\item[(U1)] There exist $E_n\uparrow E$ as $n\uparrow \infty$ and a nonnegative
  function $\fz$ such that $\sup_{i\in E_n} q_i<\infty$ and $\lim_{n\to\infty}\inf_{i\notin E_n} \fz_i=\infty$.
\item[(U2)] There exists a constant $c\in {\mathbb R}$ such that $Q\fz\le c\fz$.
\end{itemize}
}
\dethm

Certainly, for Schl\"ogl's model for instance, in condition (U2), it is more convenient to use $\ooz \fz$ instead of $Q\fz$. Besides, an important fact should
be very helpful in practice: if $\fz$ satisfies the conditions with $c\ge 0$, then so does
$M+\fz$ for every constant $M\ge 0$. In particular, a local modification of $Q$ does not interfere the conclusion.

From \rf{cmf04}{Parts I and II}, it is now clear that a large part of the theory
of $Q$-processes can be generalized to the so-called Markov jump processes on
general state space. To save the space, we will not really go to the last subject
but it is worth to mention the extension.
We now use the codes ``GS'' and ``DS'' to distinguish the ``general state space'' and the ``discrete state space'', respectively.
The sufficient part of the last theorem first appeared in \rf{cmf86a}{Theorem 2.37 (GS)} and \rf{cmf86b}{Theorem 16 (GS)}. Because it is regarded as one of the author's favourite contributions to the theory of Markov jump processes, this result was then introduced several times in the author's publications:
\rf{cmf91}{Theorem 1.11 (DS)}, \rf{cmfysj91}{Theorem 3.9 (GS)},
\rf{cmf04}{Theorem 2.25 (GS)}, \rf{cmf97}{Theorem 2.1 (DS)}, \rf{cmf05}{Theorem 9.4 (DS)}, and \rf{cmfmyh07}{Theorem 2.9 (DS)}.

 Theorem \ref{t1} is often accompanied in the publications just listed by
 the next simpler result.

 \crl\lb{t2}\;{\cms Suppose that there exist a function $\fz\ge q$ and a constant $c\in {\mathbb R}$
 such that $Q\fz\le c \fz$ on $E$. Then the $Q$-process is unique.
 }
 \decrl

 \prf Set $E_n=\{i\in E: q_i\le n\}$.
 If $M:=\sup_{i\in E}q_i<\infty$, then for large enough $n$, we have
 $E_n=E$ and so $\inf_{k\notin E_n}q_k=\infty$ by standard convention
 $\inf_{\emptyset}\fz=\infty$.
In this case, condition (U2) is trivial with $\fz=1+M$. If $M=\infty$,
then $\inf_{k\notin E_n}q_k\ge n \to \infty$ as $n\to\infty$.
Combining this with (U2), the conclusion follows from Theorem \ref{t1}.
 \deprf

Corollary \ref{t2} is almost explicit since one can simply specify
$\fz=1+ q$. This enables us to use it easier in practice.
However, such a specification makes the assumption becomes a little
stronger. We will come back this point later.

Let us make some remarks about the conditions in Theorem \ref{t1}. Condition (U2) is a relax of the
equation in (\ref{01}): finding a solution to an inequality is easier
than finding a solution to the corresponding equality. Criterion (C3) says that
there is only trivial bounded solution to the equation (\ref{01}).
Conversely, if a solution of the equation is fixed at some point, say $\uz$, such that $\fz_{\uz}=1$, then the solution
$\fz$ should be unbounded. This leads to the condition $\lim_{n\to
\infty} \inf_{k\notin E_n}\fz_k=\infty$ in (U1). Using this idea,
we prove that the assumptions in Theorem \ref{t1} are necessary for
single birth processes \rf{cmf04}{Remark 3.20}. The reason we allow
some subset of $E_n$ to be infinite is to rule out some region of $E$,
on which $\sup_{i\in E_n}q_i<\infty$. The key in the proof of this result is an economic
approximation by bounded $Q$-processes. Certainly, the necessity shows that
the assumptions of the theorem are sharp, and is valuable as illustrated
by \rf{cmf86b}{Theorem (25)}. However, it does not mean
that the inverse of the conditions in Theorem \ref{t1} can be used in practice
to show the non-uniqueness of the processes. Hence, we went to an opposite
way proving the following criterion \rf{cmf04}{Theorem 2.27 (GS), its proof in $2^{\rm nd}$
edition uses Lemma 5.18 rather than Lemma 5.15}.

\thm\lb{t4}{\cms[Non-uniqueness criterion]} {\cms For a given $Q$-matrix $Q$ on
a countable set $E$, the $Q$-processes are not unique if for some (equivalently, for all) $c>0$,
there is a bounded function $\fz$ with $\sup_{k\in E}\fz_k>0$ such that $Q\fz\ge c\fz$. Conversely, these conditions plus $\fz\ge 0$ are also necessary.
}
\dethm

We remark that three results (Theorems \ref{t1}, \ref{t4} and Corollary \ref{t2}), we have talked so far are specialized from
their original case in GS to the one in DS. Theorems \ref{t1} and \ref{t4}
are somehow the extensions of Criterion (C3)
in two opposite directions. As we will see soon that the extended
theorems are much effective than the original Criterion (C3). Using two opposite
sufficiency results instead of a single criterion is
often meaningful. For instance, for recurrence, we have a
criterion \rf{cmf04}{Proposition 4.21} which is accompanied with more practical
criteria \rf{cmf04}{Theorems 4.24 and 4.25} for the recurrence and
transience, respectively. As a companion to \rf{cmf04}{Theorem 4.25},
refer to \rf{mt93}{Theorem 8.0.2} and \rf{mh10}{Proposition 1.3} or more recent
criteria. Next, for ergodicity and nonergodicity, refer to \rf{cmf04}{Theorem 4.45\,(1)} and \rf{kl08}{Theorem 1}, respectively.
For various stability speeds/principal eigenvalues, in \ct{cmf05},
we have not only the classical variational formula, but
also dual variational formulas to describe their lower and upper bounds,
respectively.

It is interesting that there is now a direct way to prove the necessity
of Theorems \ref{t1}  in the context of DS based on a recent result by Spieksma \cite{sfm14}.

\thm\lb{t5}{\cms Everything is the same as in Theorem \ref{t1} except (U1)
is replaced by \vspace{-0.6truecm}
\begin{itemize}
\item[(U1)${}^\prime$] In the original (U1), assume in addition that each $E_n$ is finite and ignore ``$\sup_{i\in E_n} q_i<\infty$''.
\end{itemize}
}
\dethm

It is now the position to illustrate by examples the power of our results and compare conditions (U1) and (U1)${}^\prime$.

The next two examples show that in Theorem \ref{t1}, the condition ``$\lim_{n\to\infty}\fz_n$ $=\infty$'' is not necessary, which is however necessary in a criterion for recurrence used in the proof of Theorem \ref{t5} (see its proof below).

\xmp\lb{t6} {\cms Let $E$ be a countable set and $Q=(q_{ij})$ be a bounded conservative
$Q$-matrix on $E$. Then assumptions of Theorem \ref{t1} hold but its test function $\fz$ can be bounded.}
\dexmp

\prf (a) Simply set $E_n\equiv E$ (may be infinite) for every $n\ge 1$ and $\fz_i\equiv 1$. Then it is obvious
that $0=Q\fz \le \fz$ and $\lim_n \inf_{i\notin E_n}\fz_i=\infty$ since $\inf_{\emptyset}\fz=\infty$
by the standard convention. Hence by Theorem \ref{t1}, the process is unique. As we have seen before, Corollary \ref{t2} is also applicable in such a trivial case.

(b) Knowing that the process is unique, then by Theorem \ref{t5}, there should exist a $\fz$ satisfying (U1)${}^\prime$, as well as (U2). The problem is that the resulting $\fz$ is not explicitly known when $E$ is infinite. In this sense, Theorem \ref{t5} is theoretic correct but not practical in such simplest case.\deprf

\xmp\lb{t7} {\cms Let $E={\mathbb Z}_+$ and $Q^{(1)}$ be a bounded conservative
$Q$-matrix on $E$. Denote its test function by $\fz^{(1)}\equiv 1$
as in the last example.
Next, let $Q^{(2)}$ be a conservative $Q$-matrix on $E$ satisfying the assumptions of Theorem \ref{t1} with a sequence of finite subsets $\{E_n\}_{n\ge 1}$
and a test function $\fz^{(2)}$.
Finally, we construct a new $Q$ as follows: on the odd numbers in $E$, we
use the transition mechanism of $Q^{(1)}$, and on the even numbers in $E$,
we adopt the one of $Q^{(2)}$. Define
$\fz=\fz^{(1)}$ on the odd numbers and $\fz=\fz^{(2)}$ on the even numbers.
Then the assumptions of Theorem \ref{t1} hold but its test function $\fz_n$ has no limit
as $n\to\infty$: $\varlimsup_{n\to\infty}\fz_n=\infty$ and $\varliminf_{n\to\infty}\fz_n=1$.
}\dexmp

\prf First, note that for the original $Q^{(2)}$ on $E$, because each $E_n$
is a finite subset of $E$, the condition $\lim_{n\to\infty}\inf_{k\notin E_n}\fz_k^{(2)}=\infty$
is equivalent to $\lim_{n\to\infty}\fz_n^{(2)}=\infty$. Therefore, we have
$$\varlimsup_{n\to\infty}\fz_n=\lim_{n\to\infty}\fz_n^{(2)}=\infty,\qquad
\varliminf_{n\to\infty}\fz_n=\lim_{n\to\infty}\fz_n^{(1)}=1.$$
To show the assumptions in Theorem \ref{t1} hold, simply let
$E_0=\{\text{odd inte\-gers}\}$, and let $E_n\,(n\ge 1)$ be the union of
$E_0$ and the natural modification of the original $E_n$ used for $Q^{(2)}$.
Then the resulting $E_n\uparrow E$ as $n\to\infty$,
$\sup_{k\in E_n}q_k<\infty$ for each $n\ge 0$,
and $$\lim_{n\to\infty}\inf_{k\notin E_n}\fz_k=
\lim_{n\to\infty}\inf_{k\notin E_n}\fz_k^{(2)}
=\lim_{n\to\infty}\fz_n^{(2)} =\infty.$$
Finally, because of the independence
of $Q^{(1)}$ and $Q^{(2)}$,  $\fz^{(1)}$ and $\fz^{(2)}$,
the condition $Q\fz \le \max\{c_2,\, 1\}\fz$ on the set of odd numbers follows from
$$Q^{(1)}\fz^{(1)} \le  \fz^{(1)}
\quad \text{on } \;  E;$$
and the same condition on the set of even numbers follows from
$$Q^{(2)}\fz^{(2)} \le c_2 \fz^{(2)}\quad \text{on } \;  E.$$
We have thus obtained the required conclusion.

As mentioned in the last proof, in the present situation, we do not
know how to use Theorem \ref{t5}.\deprf

Note that the last matrix $Q$ is reducible. However, we can add a connection between 0 and 1 to produce an irreducible version of the example. This is not essential since a local modification
does not interfere the uniqueness problem. Furthermore, one may replace
the set $\{\text{odd integers}\}$ or $\{\text{even integers}\}$ by any infinite subset of $E$, but not $E$ itself, the set of primer numbers for instance.
The conclusion of Example \ref{t7} remains the same by an obvious modification.

The point is that some $E_n$ is allowed to be infinite in (U1) but not in (U1)${}^\prime$.

\xmp\lb{t8}{\cms The pure birth process is unique iff (\ref{02}) holds. In particular, set $q_{n, n+1}=$\,the $n$th primer, then Theorem \ref{t1} is
suitable but Corollary \ref{t2} fails.}
\dexmp

\prf Note that $q_k=q_{k, k+1}$ for $k\ge 0$.

(a) If $\sum_k q_k^{-1}=\infty$, set $E_n=\{0, 1, \ldots, n\}$ and
$$\fz_k=1+\sum_{1\le j\le k-1}\frac{1}{q_j}\to \infty \quad \text{as }k\to\infty.$$
Then $Q\fz\le \fz$ and so Theorem \ref{t1} gives us the uniqueness of the
processes. In the particular case that $q_{n, n+1}=n+1$, the above $\fz$
has order $\log n$. However, we can also choose $\fz_n=1+n$ and apply Theorem \ref{t1}. This shows that there are some freedom in choosing
$\fz$.

(b) If $M:=\sum_k q_k^{-1}<\infty$, set $E_n$ as above and
$$\fz_k=
\frac 1 2+\sum_{1\le j\le k-1}\frac{1}{q_j}-M\in \bigg[\frac 1 2 -M,\; \frac 1 2\bigg].$$
Then $\sup_k \fz_k=1/2>0$, $Q\fz \ge \fz$, and so by Theorem \ref{t4},
the processes are not unique. We remark that it would be awful to use the necessity in Theorems \ref{t1} or \ref{t5} to prove this non-uniqueness property.

(c) The last assertion is due to J.L. Zheng (cf. \rf{cmf86a}{Example 2.3.12} or \rf{cmf04}{Example 2.26}).\deprf

\medskip

\nnd{\bf Proof of the uniqueness for Example \ref{t0}}.
For $i\in E={\mathbb Z}_+^S$, define its level by $|i|=\sum_{u\in S}i_u$
and set $E_n=\{i\in E: |i|\le n\}$ for $n\ge 1$.

(a) Next, define
$\fz(i)=1+|i|$. Then it is clear that $\lim_{n\to\infty}\inf_{k\notin E_n}\fz(k)=\infty$. Because the diffusions do not change the levels,
we have
$$\ooz \fz(|i|)=\sum_{u\in S} [b(i_u)-a(i_u)]=\sum_{u\in S}\big[\az_0-\az_1i_u+\az_2 i_u^2-\az_3 i_u^3\big]
$$
for some positive $\{\az_k\}_{k=0}^3$. Next, since
$$\sum_{u\in S}i_u^2\le |i|^2, \qqd
\frac{1}{|S|}\sum_{u\in S}i_u^3\ge \bigg(\frac{|i|}{|S|}\bigg)^3
\;\text{(Jensen's inequality)},$$
where $|S|$ is the candinality of $S$ (finite but arbitrary), we have
$$\ooz \fz(|i|)\le \az_0'-\az_1' |i|+\az_2'|i|^2-\az_3'|i|^3$$
for some positive $\{\az_k'\}_{k=0}^3$.
Now, because the right-hand side becomes negative for large enough
$|i|$, it is clear that
$\ooz \fz(|i|)\le c\fz(|i|)$ for every $i\in E$ and large enough $c$.
The assertion now follows from Theorem \ref{t1}.
Hopefully, we have seen the role played by the geometry of $E$.
The proof shows the power of our result. A good sufficiency result
may be more effective than a criterion.

(b) It is also possible to use Corollary \ref{t2} to prove the required assertion, simply choose
$\fz (i)=\gz\big(1+ \sum_{u\in S}i_u^3\big)$.
First, choose $\gz$ large enough so that $\fz\ge q$. Next, choose $c$
large enough so that $\ooz \fz\le c\fz$.
\qed
\smallskip

It is worthy to mention that in accompany to Theorem \ref{t1}, we also have a similar,
practical sufficiency result for (exponential) ergodicity. Refer to \rf{cmf89}{Theorem 3 (GS)},
\rf{cmf91}{Theorem (1.18) (DS)}, \rf{cmf04}{Corollary 4.49 (DS) and Theorem 14.1 (GS)}.

In the past nearly 30 years, Theorem \ref{t1} and Corollary \ref{t2}
have very successful applications. A list of the literature was collected in
\rf{cmf05}{\S 9.2}. Certainly, the results used a lot by the author
(in \cite{cmf04} for instance). In particular, it was used at the first step
to construct a large class of infinite-dimensional processes (\rf{cmf04}{\S 13.2}), 15 models are included
in \rf{cmf04}{\S 13.4}. Corollary \ref{t2} with some extension was used by Song (1988) \cite{sjs88} in a quite earlier stage for Markov decision processes
moving from bounded to unbounded situation.
It is now quite often to see the influence of the study on Markov Jump
processes to the theory of Markov decision processes.
Based on \cite{cmf86b},
Theorem \ref{t1} was collected into Anderson \rf{awj91}{Corollary 2.2.16},
its originality was unfortunately ignored, even though the original paper
\ct{cmf86b} is included in the references of the book.
For some corrections and comments on the last book, refer to \cite{cmf96}.
Very recently, Theorem \ref{t1}(GS) is applied by Chen and Ma \cite{chxmzm14}
to genetic study having continuous state space. Finally, we mention that the
results have already
extended to the time-inhomogeneous case by Zheng \cite{zjlzxg87} and \cite{zjl93}
using the martingale approach.

Before going to the proofs, note that equation (\ref{01}) is equivalent to
\be \ppz (\lambda) u=u, \quad 0\le u\le 1\;\;\text{on }E,\qquad \lz>0, \lb{03}\de
where
$$\ppz (\lambda)=\bigg(\frac{(1-\dz_{ij})q_{ij}}{\lambda+q_i}: \; i, j\in E\bigg).$$
Here the matrix $\ppz (\lambda)$ is sub-stochastic. We introduce a fictitious state $\ddz$ and define on the enlarged state space $E_{\ddz}=E\cup \{\ddz\}$ a new transition probability matrix
$$\ppz_{ij}^{\ddz} (\lambda)
=\begin{cases}
\ppz_{ij} (\lambda)  \quad &\text{if } i, j\in E\\
\frac{\lz}{\lz+q_i} &\text{if } i\in E, \; j=\ddz\\
p_j &\text{if } i=\ddz, j\in E
\end{cases}$$
where $(p_j: j\in E)$ is a positive probability measure on $E$. The enlarged
transition probability matrix is irreducible even the original one may be not.

\lmm\lb{t9}{\cms The equation (\ref{01}) has zero solution only iff so does the equation
\be \ppz^{\ddz} (\lambda) \big(u{\mathbbold{1}}_E\big)=u, \quad 0\le u\le 1\;\;\text{\cms on }E_{\ddz},\qquad  \lz>0.\lb{04}\de
Thus, the original $Q$-process is unique iff the $\ppz^{\ddz} (\lambda)$-chain is recurrent.
}
\delmm

\prf Noting that $u_{\ddz}=\sum_{k\in E}p_k u_k$, it is clear that
$u_{\ddz}=0$ iff $u_k=0$ for all $k\in E$ since $p_k>0$ for all $k\in E$.
Equation (\ref{04}) restricted to $E$ coincides with (\ref{03}) and then
(\ref{01}). This proves the first assertion.

To prove the second assertion, it suffices to note that $\ppz^{\ddz} (\lambda)$-chain is recurrent iff equation (\ref{04}) has only trivial solution. The last result comes from \cite{ysjcmf86}, \rf{cmf86a}{Lemma 12.1.27}, or \rf{cmf04}{Lemma 4.51}. We remark here that the regularity assumption
used in the cited references can be replaced by the minimal process, due to
the equivalence of recurrence of the minimal process and its embedded chain.
Refer to \rf{cmf86a}{Lemma 12.3.1}, or \rf{cmf04}{Theorem 4.34}.
\deprf


\nnd{\bf Proof of Theorem \ref{t5}}. When $|E|<\infty$, the conclusion is trivial and the assumptions hold for the specific $E_n\equiv E$ and $\fz_i\equiv 1$ as seen from proof (a) of Example \ref{t6}. Hence we may assume that $E={\mathbb Z}_+$. Since each $E_n$ is finite, the
condition $\lim_n \inf_{k\notin E_n}\fz_k=\infty$ becomes
$\lim_{n\to\infty} \fz_n=\infty$. In this case, conditions (U1)${}^\prime$ and (U2)
consist a criterion for the recurrence of the Markov chain
$\ppz^{\ddz} (\lambda) $, refer to \rf{cmf04}{Theorem 4.24} and its references within.

We remark that it is at this point, the finiteness of $E_n$ is required and so the present sufficiency proof is not suitable for Theorem \ref{t1}. At the moment, we do not know how to extend the necessity result of Theorem \ref{t5} from DS to GS.

Here is a part of an alternative proof given in \cite{sfm14}.
Let $P^{\min}(\lz)$ be the Laplace transform of $P^{\min}(t)$. Using the
second successive approximation scheme for the backward Kolmogorov equation
(goes back to \rf{wf40}{Theorem 1}), we obtain
$$P^{\min}(\lz)= \sum_{n=0}^{\infty}\ppz(\lz)^n \text{diag}\bigg(\frac{1}{\lz+q}\bigg)$$
(cf. \rf{cmf04}{page 75, line -6}). Hence
$$\lz P^{\min}(\lz)\,\text{column}\,(1)=\sum_{n=0}^{\infty}\ppz(\lz)^n \, \text{column}\bigg(\frac{\lz}{\lz+q}\bigg).$$
The process is unique iff the left-hand side equals 1 at some/every $i\in E$,
the right-hand side is the probabilistic decomposition of the
time that the Markov chain $\ppz^{\ddz} (\lambda)$ starts from some $i\in E$,
first visits $\ddz$ at some step $n\ge 1$, which
equals 1 iff the irreducible Markov chain $\ppz^{\ddz} (\lambda)$ is recurrent.
We have thus come back to the last lemma.\qed
\medskip

\nnd{\bf Proof of Theorem \ref{t1}}. Here we adopt a circle argument.

(U1)${}^\prime$$\,+\,$(U2)\,$\Longrightarrow$\,(U1)$\,+\,$(U2). This is easy since
(U1) is weaker that (U1)${}^\prime$.

(U1)$\,+\,$(U2)\,$\Longrightarrow$\, uniqueness. This is the sufficiency
 part of Theorem \ref{t1} and was proved long time ago, even for GS.

 Uniqueness\,$\Longrightarrow$\,(U1)${}^\prime$$\,+\,$(U2). This is the necessity
 part of Theorem \ref{t5}.
\qed

\medskip

We remark that a similar phenomena is appeared in Theorem \ref{t4}, the conditions for sufficiency are weaker than the ones for necessity. As
we have seen from Example \ref{t8}, this is very helpful in practice. However, these
conditions are actually equivalent: conditions for necessity $\Longrightarrow$ conditions for sufficiency $\Longrightarrow$ non-uniqueness $\Longrightarrow$ conditions for necessity.

In view of these discussions, one may combine Theorems \ref{t1} and \ref{t5} into one having the style of Theorem \ref{t4}.

In conclusion, this note as well as the practice during the past 30 years
confirm that the sufficient part of Theorem \ref{t1} and Theorem (Criterion) \ref{t4} are not only powerful but also sharp, even though at the moment
we are still unable to prove the necessity part of Theorem \ref{t1} for
general state spaces.

\medskip

\nnd{\bf Acknowledgments}. {\small
The author thanks Yong-Hua Mao for bringing \cite{sfm14} to the attention.
Research supported in part by the
         National Natural Science Foundation of China (No. 11131003),
         the ``985'' project from the Ministry of Education in China,
and the Project Funded by the Priority Academic Program Development of
Jiangsu Higher Education Institutions.
}

\medskip

\nnd {\small
Mu-Fa Chen\\
School of Mathematical Sciences, Beijing Normal University,
Laboratory of Mathematics and Complex Systems (Beijing Normal University),
Ministry of Education, Beijing 100875,
    The People's Republic of China.\newline E-mail: mfchen@bnu.edu.cn\newline Home page:
    http://math.bnu.edu.cn/\~{}chenmf/main$\_$eng.htm}


\begin{thebibliography}{10}
\setlength{\itemsep}{-0.8ex}
{\small

\bibitem{awj91}
Anderson, W.J. (1991).
{\it Continuous-Time {M}arkov Chains}.
Springer Series in Statistics. \lb{awj91}

\bibitem{cmf85}
Chen, M.F. (1985).
{\it Infinite-dimensional reaction-diffusion processes}.
Acta Math. Sin. (N.S.) 1(3), 261--273.  \lb{cmf85}

\bibitem{cmf86a}
Chen, M.F. (1986a).
{\it Jump Processes and Interacting Particle Systems {\rm(in
  Chinese)}}.
Beijing Normal Univ. Press. \lb{cmf86a}

\bibitem{cmf86b}
Chen, M.F. (1986b).
{\it Couplings of jump processes}.
Acta Math. Sinica, New Series 2(2), 123--136. \lb{cmf86b}

\bibitem{cmf89}
Chen, M.F. (1989).
{\it Stationary distributions of infinite particle systems with non-compact
state spaces}.
Acta Math. Sci. 9(1), 7--19.\lb{cmf89}

\bibitem{cmf91}
Chen, M.F. (1991).
{\it On three classical problem for {M}arkov chains with continuous time
  parameters}.
J. Appl. Prob. 28(2), 305--320. \lb{cmf91}

\bibitem{cmf96}
Chen, M.F. (1996).
{\it A comment on the book ``Continuous-Time Markov Chains'' by W.J. Anderson}.
Chin. J. Appl. Prob. Stat. 12(1), 55--59.\newline
See also  arXiv:1412.5856.  \lb{cmf96}

\bibitem{cmf97}
Chen, M.F. (1997).
{\it Reaction-diffusion processes}.
    Chin. Sci. Bull. 42(23), 2465--2474 $(${Chin. Ed.}$)$; 43(17), 1998, 1409--1421  $(${\rm Eng. Ed.}$)$. \lb{cmf97}

\bibitem{cmf04}
Chen, M.F. (2004).
{\it From Markov Chains to
Non-equilibrium Particle Systems}.
World Scientific. 2$^{\text{nd}}$ ed. (1$^{\text{st}}$ ed., 1992). \lb{cmf04}

\bibitem{cmf05}
Chen, M.F. (2005).
{\it Eigenvalues, Inequalities, and Ergodic Theory}.
Springer, London.  \label{cmf05}

\bibitem{cmfmyh07}
Chen. M.F. and Mao, Y.H. (2007).
{\it Introduction to Stochastic Processes {\rm(in
  Chinese)}}.
Heigher Edu. Press. \lb{cmfmyh07}

\bibitem{cmfysj91}
Chen. M.F. and Yan, S.J. (1991).
{\it Jump processes and particle systems},
in ``Probability Theory and its Applications in China'', edited by S.J.
      Yan, C.C. Yang and J.G. Wang, Providence,  AMS. 118, 23-57. \lb{cmfysj91}

\bibitem{chzh}
Chen, M.F. and Zhang, Y.H. (2014).
{\it Unified representation of formulas for single birth processes}.
Front. Math. China 9(4), 761--796. \lb{chzh}

\bibitem{chxmzm14}
Chen, X. and Ma, Z.M. (2014).
{\it A transformation of Markov jump processes and applications in genetic study}.
Discrete and continuous dynamical systems 34(12), 5061--5084.  \lb{{chxmzm14}}

\bibitem{ckl67}
Chung, K.L.(1967).
{\it Markov Chains with Stationary Transition Probabilities}.
1$^{\rm st}$ ed., 1960, 2$^{\rm nd}$ ed., 1967, Springer-Verlag, New York.\lb{ckl67}

\bibitem{drl52}
Dobrushin, R.L. (1952).
{\it On conditions of regularity of stationary Markov processes with a denumerable number of possible states \text{\rm (in Russian)}}.
Uspehi Matem. Nauk (N.S.) 7(6), 185--191. See http://www.cpt.univ-mrs.fr/dobrushin/list.html \lb{drl52}
{\vspace{-0.45truecm}}
\bibitem{wf40}
Feller, W. (1940).
{\it On the integro-differential equations of pure discontinuous Markov Processes}
Trans. Amer. Math. Soc. 48, 488--515. \lb{wf40}

\bibitem{wf57}
Feller, W. (1957).
{\it  On boundaries and lateral conditions for Kolmogorof differential
equations}.
Ann. Math. 65, 527--570. \lb{wf57}

\bibitem{gs75}
Gikhman, I.I., Skorokhod, A.V. (1975).
{\it The Theory of Stochastic Processes }II.
Springer, New York. \lb{gs75}

\bibitem{mh10}
Hairer, M. (2010).
{\it Convergence of Markov Processes}. Lecture Notes.\newline
http://www.hairer.org/notes/Convergence.pdf  \lb{mh10}

\bibitem{kl08}
Kim, B. and Lee, I. (2008).
{\it Tests for nonergodicity of denumerable continuous time Markov processes}.
Comput. Math. Appl. 55: 1310--1321.  \lb{kl08}

\bibitem{mt93}
Meyn, S.P. and Tweedie, R.L. (2009).
{\it Markov Chains and Stochastic Stability \text{\rm 2$^{\rm nd}$ ed}}.
Cambridge Univ. Press. \lb{mt93}

\bibitem{gehr57}
Reuter, G.E.H. (1957)
{\it Denumerable Markov Processes}.
Acta Math. 57, 1--46. \lb{gehr57}

\bibitem{sjs88}
Song, J.S. (1988).
{\it Continuous time Markov decision processes with non-uniformly
  bounded transition rates}.
Sci. Sin. Ser. A. 12(11), 1281--1291.  \lb{sjs88}

\bibitem{sfm14}
Spieksma, F.M. (2014).
{\it Countable state Markov processes:
non-explosiveness and moment function}.
To appear in Probab. Eng. and Inform. Sci.\lb{sfm14}

\bibitem{ysjcmf86}
Yan, S.J. and Chen, M.F. (1986)
{\it multidimensional $Q$-processses}.
Chin. Ann. Math. 7B(1), 90--110. \lb{ysjcmf86}

\bibitem{zjl93}
Zheng, J.L.(1993).
{\it Phase Transitions of Ising Model on Lattice Fractals,
Martingale Approach for $q$-processes  \rm{(in Chinese)}}.
Ph.D. Thesis, Beijing Normal Univ. \lb{zjl93}
{\vspace{-0.45truecm}}
\bibitem{zjlzxg87}
Zheng, J.L. and Zheng, X.G.(1987).
{\it A martingale approach for $Q$-processes {\rm (abstract)}}.
Sci. Bull. 21, 1457--1459.  \lb{zjlzxg87}
}
\end{thebibliography}
\end{document}